\documentclass[12pt]{elsarticle}

\usepackage{graphicx}
\usepackage{amsfonts}
\usepackage{amssymb}
\usepackage{amsmath}
\usepackage{epic}
\usepackage{eepic}
\usepackage{srcltx} 
\usepackage{dsfont}
\usepackage{color}
\usepackage{stmaryrd}
\usepackage{hhline}

\usepackage{amsthm}

\begin{document}

\begin{frontmatter}



\title{Improved accuracy for time-splitting methods for the numerical solution of parabolic equations}

\author[a]{A. Arrar\'as\corref{cor1}}
\ead{andres.arraras@unavarra.es}
\author[a]{L. Portero}
\ead{laura.portero@unavarra.es}
\cortext[cor1]{Corresponding author.}
\address[a]{Departamento de Ingenier\'{\i}a Matem\'atica e Inform\'atica, Universidad P\'ublica de Navarra, Edificio Las Encinas, Campus de Arrosad\'{\i}a, 31006 Pamplona (Spain)}

\begin{abstract}
In this work, we study time-splitting strategies for the numerical approximation of evolutionary reaction--diffusion problems. In particular, we formulate a family of domain decomposition splitting methods that overcomes some typical limitations of classical alternating direction implicit (ADI) schemes. The splitting error associated with such methods is observed to be $\mathcal{O}(\tau^2)$ in the time step $\tau$. In order to decrease the size of this splitting error to $\mathcal{O}(\tau^3)$, we add a correction term to the right-hand side of the original formulation. This procedure is based on the improved initialization technique proposed by Douglas and Kim in the framework of ADI methods. The resulting non-iterative schemes reduce the global system to a collection of uncoupled subdomain problems that can be solved in parallel. Computational results comparing the newly derived algorithms with the Crank--Nicolson scheme and certain ADI methods are presented.
\end{abstract}

\begin{keyword}
alternating direction implicit \sep domain decomposition \sep partition of unity \sep splitting error \sep time-splitting method

\vspace*{0.25cm}
\MSC[2010] 35K57 \sep 65M55 \sep 68W10

\end{keyword}

\end{frontmatter}



\section{Introduction}\label{sec:introduction}

The numerical approximation of parabolic problems using time-splitting procedures has been a wide field of research since the pioneering works of Douglas, Peaceman and Rachford in the decade of the 1950s (cf. \cite{dou:55,dou:pea:55,dou:rac:56,pea:rac:55}). In such works, they introduced the so-called alternating direction implicit (ADI) methods by noting that, in a $d$-dimensional spatial domain $\Omega$, the diffusion operator $-\nabla\cdot(a\nabla)$ can be expressed as the sum of $d$ one-dimensional operators $\{\partial_k(a\partial_k)\}_{k=1,2,\ldots,d}$, $a$ being a uniformly positive function on $\overline{\Omega}$ and $\partial_k=\partial/\partial x_k$. Using this idea, multidimensional parabolic problems can be solved as a sequence of one-dimensional problems, each formulated for one of the spatial variables under consideration. Any time-stepping procedure based on a component-wise splitting of this kind is called a locally one-dimensional method (see, e.g., \cite{hun:ver:03}). For an extensive study of ADI and related time-splitting methods, we refer to the monographs \cite{far:hav:09,hun:ver:03,mar:90}.

Significantly, an ADI method can be viewed as a perturbation of some underlying implicit scheme, such as the backward Euler or the Crank--Nicolson method. In this setting, the ADI method may be formulated as the corresponding implicit scheme plus a perturbation term called the splitting error. In general, this splitting error is of the same --or higher-- order in the time step $\tau$ as the truncation error associated with the underlying unsplit method. As a result, the asymptotic rate of convergence for both the ADI and its underlying method should be the same; in practice, however, the actual errors associated with the former are much larger than those associated with the latter. This fact is due to the presence of the splitting error and, typically, it is considered to be the main drawback of time-splitting methods. In order to reduce the size of such an error, Douglas and Kim introduced in \cite{dou:kim:01} the so-called alternating direction method with improved initialization (AD-II) (sometimes referred to as the modified alternating direction (AD-M) method, see \cite{arb:hua:yan:07}). Essentially, they proposed to add a correction term to the right-hand side of the original ADI scheme with the aim of reducing the splitting error from $\mathcal{O}(\tau^2)$ to $\mathcal{O}(\tau^3)$. This idea was used in \cite{arb:hua:yan:07} to formulate improved ADI methods for regular and mixed finite element discretizations, and further studied in \cite{kar:09} to derive linear multistep methods by the approximate factorization technique.

In this paper, we extend the improved initialization procedure of Douglas and Kim to the case of domain decomposition-based splittings. This kind of splittings was first introduced in \cite{vab:89,vab:94} to obtain the so-called regionally-additive schemes, and has been subsequently studied in \cite{mat:pol:rus:wan:98,por:buj:jor:04,vab:08} for the solution of evolutionary problems. In the context of linear and semilinear parabolic equations, it has been successfully used in combination with various spatial discretization techniques, such as mimetic finite differences (cf. \cite{arr:por:jor:09,por:arr:jor:11}), mixed finite elements (cf. \cite{arr:por:14}) or multipoint flux approximation methods (cf. \cite{arr:por:yot:14}). Some additional results regarding nonlinear parabolic equations can be found in \cite{arr:por:jor:10,por:arr:jor:10}. The monographs \cite{mat:08,sam:mat:vab:02} show an overview of some recent contributions to the topic.

The key to the efficiency of domain decomposition splitting methods lies in reducing the system matrix to a collection of uncoupled submatrices of lower dimension. As compared to classical overlapping domain decomposition algorithms (cf. \cite{qua:val:99}), this approach does not involve any Schwarz iteration procedure, thus reducing the computational cost of the overall solution process. In addition, it overcomes two typical limitations of alternating direction splittings, namely: (a) their need to deal with rectangular or hexahedral spatial grids (in two- or three-dimensional problems, respectively); and (b) their difficulty to handle mixed derivative terms. In this respect, although several ADI methods have been specifically designed along the years to overcome this latter constraint (see, e.g., \cite{dou:gun:64,mck:mit:70,mck:wal:wil:96} or, more recently, \cite{int:mis:13,int:wel:07,int:wel:09}), no AD-II scheme with this property has been developed so far.

In order to introduce the improved time-splitting procedures, we consider a parabolic initial-boundary value problem of the form
\begin{subequations}\label{continuous:problem}
\begin{align}
u_t-\nabla\cdot(a\nabla u)+cu&=f,&&\hspace*{-1.5cm}\mathbf{x}\in\Omega, &&\hspace*{-1.5cm}0<t\le T,\label{continuous:problem:a}\\
u&=0,&&\hspace*{-1.5cm}\mathbf{x}\in\partial\Omega, &&\hspace*{-1.5cm}0<t\le T,\label{continuous:problem:b}\\
u&=u_0,&&\hspace*{-1.5cm}\mathbf{x}\in\Omega, &&\hspace*{-1.5cm}t=0,\label{continuous:problem:c}
\end{align}
\end{subequations}
where $\Omega\subset\mathds{R}^2$ is a bounded open domain with boundary $\partial\Omega$, $a(\mathbf{x})\in\mathds{R}^{2\times2}$ is a symmetric positive definite matrix function, with elements $\{a_{ij}(\mathbf{x})\}_{i,j=1,2}$, $c(\mathbf{x})$ is a uniformly positive function on $\overline{\Omega}$, and the subscript $t$ denotes partial differentiation with respect to time. The entries of $a(\mathbf{x})$ and the terms $c(\mathbf{x})$, $f(\mathbf{x},t)$ and $u_0(\mathbf{x})$ are assumed to be sufficiently smooth. In the sequel, we denote by $Au=-\nabla\cdot(a\nabla u)+cu$ the elliptic operator applied to the exact solution $u(\mathbf{x},t)$. For simplicity in the exposition, we consider homogeneous Dirichlet boundary conditions, although more general conditions can also be handled.

For the sake of completeness, it must be mentioned that the solution of problem (\ref{continuous:problem}) can also be approximated with a high order of accuracy by means of unsplit methods. Among the time-stepping schemes that have been recently designed for this purpose, it is worth referring to the so-called ADER (Arbitrary DERivative in space and time) approach. This method was initially developed to provide high-order approximations to the solution of hyperbolic conservation laws (see, e.g., \cite{dum:ena:tor:08,tit:tor:02}), and was subsequently formulated for nonlinear reaction--diffusion problems (cf. \cite{tor:hid:09}). More recently, it has been extended in \cite{hid:dum:11} to nonlinear systems of advection--diffusion--reaction equations involving stiff source terms. Some novel contributions to the topic can be found in \cite{mon:tor:14,tor:mon:14}.

The rest of the paper is organized as follows. In \S\ref{sec:unsplit}, we formulate some classical unsplit implicit methods for the numerical solution of problem (\ref{continuous:problem}). Two time-splitting strategies are described in \S\ref{sec:split}: on one hand, the well-known Douglas (cf. \cite{dou:62}) and Douglas--Rachford (cf. \cite{dou:rac:56}) alternating direction methods; on the other, a non-iterative domain decomposition splitting technique based on a family of partition of unity functions. The splitting error associated with the preceding schemes is studied in \S\ref{sec:splitting:error}. In this section, we also introduce a correction procedure to reduce the size of such an error. Numerical experiments comparing the previous time-splitting schemes with the Crank--Nicolson method are reported in \S\ref{sec:experiments}. The paper ends with some concluding remarks summarized in \S\ref{sec:conclusions}.

\section{Classical implicit methods}\label{sec:unsplit}

Let us consider a suitable partition $\mathcal{T}_h$ of the spatial domain $\Omega$, where $h$ denotes the maximal grid spacing. In principle, $\mathcal{T}_h$ is an unstructured grid composed of either triangular or quadrilateral elements. We will see in the next section that, if an ADI method is used for solving (\ref{continuous:problem}), $\mathcal{T}_h$ must be assumed to be a rectangular grid. Domain decomposition splittings, by contrast, will remain valid for unstructured partitions. On the other hand, let us consider a constant time step $\tau>0$ and define the discrete times $t_n = n\tau$, for $n=0,1,\ldots,N_T+1$, with $N_T=[T/\tau]-1$.

In this framework, let $A_h$ be a symmetric positive definite matrix obtained from finite difference or finite element discretization of the elliptic operator $A$ with order of accuracy $\mathcal{O}(h^s)$. Then, if $U^n_h$ denotes the fully discrete solution at time $t_n$, the classical implicit time-stepping schemes can be written together, for $n=0,1,\ldots,N_T$, in the form
\begin{equation}\label{be:cn}
\frac{U_h^{n+1}-U_h^n}{\tau}+A_h(\theta U_h^{n+1}+(1-\theta)U_h^n)=F_h^{n+\theta},
\end{equation}
where $\theta=1$ for backward Euler and $\theta=\frac{1}{2}$ for Crank--Nicolson. Note that $F_h^{n+\theta}=\theta F_h(t_{n+1})+(1-\theta)F_h(t_{n})$, $F_h(t_n)$ being the discrete forcing term at time $t_n$. The initial condition is given by $U_h^0=\mathcal{P}_hu_0$, where $\mathcal{P}_h$ is a suitable restriction or projection operator. It is well known that, for a spatial discretization of order $\mathcal{O}(h^s)$, the backward Euler solution converges to the true solution with order $\mathcal{O}(\tau+h^s)$; in turn, the Crank--Nicolson solution converges with order $\mathcal{O}(\tau^2+h^s)$. In compact form, both results can be written together as $\mathcal{O}(\tau^{3-2\theta}+h^s)$.

\section{Time-splitting methods}\label{sec:split}

In this section, we introduce two classes of time-splitting procedures for the solution of problem (\ref{continuous:problem}), namely: (a) the classical alternating direction methods based on a component-wise splitting; and (b) a family of non-iterative schemes based on an overlapping domain decomposition splitting.

\subsection{Alternating direction splitting}

Let us assume that $a(\mathbf{x})$ in (\ref{continuous:problem:a}) is a diagonal matrix with elements $a_{11}(\mathbf{x})$ and $a_{22}(\mathbf{x})$. This assumption avoids the existence of mixed derivative terms in the parabolic problem. In this way, the elliptic operator $A$ can be expressed as the sum $A=A_1+A_2$, where
\begin{equation}\label{adi:splitting}
A_{1}u=-(a_{11}u_x)_x+\dfrac{1}{2}\,cu,\qquad
A_{2}u=-(a_{22}u_y)_y+\dfrac{1}{2}\,cu.
\end{equation}
If the spatial domain $\Omega$ admits a rectangular grid $\mathcal{T}_h$, problem (\ref{continuous:problem}) can be approximated, for $n=0,1,\ldots,N_T$, by an ADI method of the form
\begin{equation}\label{adi:scheme}
\begin{aligned}
\frac{V_h^{n,1}-V_h^n}{\tau}+A_{1h}(\theta V_h^{n,1}+(1-\theta)V_h^n)+A_{2h}V_h^n&=F_h^{n+\theta},\\
\frac{V_h^{n+1}-V_h^n}{\tau}+A_{1h}(\theta V_h^{n,1}+(1-\theta)V_h^n)+A_{2h}(\theta V_h^{n+1}+(1-\theta)V_h^n)&=F_h^{n+\theta}.
\end{aligned}
\end{equation}
Here, $A_{1h}$ and $A_{2h}$ are suitable discretizations of $A_{1}$ and $A_{2}$, respectively, such that $A_h=A_{1h}+A_{2h}$. Note that we introduce the notation $V_h^n$ to distinguish the ADI discrete solution from the unsplit discrete solution $U_h^n$ of the preceding section. In this case, the initial condition is also given by $V_h^0=\mathcal{P}_hu_0$.

The previous discretization scheme yields some classical ADI methods for certain values of the parameter $\theta$. More precisely, if $\theta=1$, we recover the Douglas--Rachford method, first proposed in \cite{dou:rac:56}; i.e., for $n=0,1,\ldots,N_T$,
\begin{equation}\label{douglas:rachford}
\begin{aligned}
(I+\tau A_{1h})\,V_h^{n,1}&=(I-\tau A_{2h})\,V_h^n+\tau F_h^{n+1},\\
(I+\tau A_{2h})\,V_h^{n+1}&=V_h^{n,1}+\tau A_{2h}V_h^n,
\end{aligned}
\end{equation}
for a given $V_h^0$, where $I$ denotes the identity matrix. In turn, if $\theta=\frac{1}{2}$, the Douglas method (cf. \cite{dou:62}) is obtained; i.e., for $n=0,1,\ldots,N_T$,
\begin{equation}\label{douglas:gunn}
\begin{aligned}
\left(I+\dfrac{\tau}{2} A_{1h}\right)V_h^{n,1}&=\left(I-\dfrac{\tau}{2} A_{1h}-\tau A_{2h}\right)V_h^n+\tau F_h^{n+\frac{1}{2}},\\
\left(I+\dfrac{\tau}{2} A_{2h}\right)V_h^{n+1}&=\textstyle V_h^{n,1}+\dfrac{\tau}{2} A_{2h}V_h^n,
\end{aligned}
\end{equation}
for a given $V_h^0$.

\subsection{Domain decomposition splitting}

As an alternative to component-wise splitting methods, we present a time-stepping procedure based on a suitable decomposition of the spatial domain. Unlike the ADI methods introduced above, this technique permits to handle mixed derivative terms (i.e., full tensor coefficients $a(\mathbf{x})$ in (\ref{continuous:problem:a})), and is also valid for unstructured partitions $\mathcal{T}_h$ of the spatial domain $\Omega$. In this case, we consider the splitting $A=A_{1}+A_{2}+\ldots+A_{m}$ into an arbitrary number $m$ of split terms, where
\begin{equation}\label{dd:splitting}
A_{k}u=\textstyle -\nabla\cdot(\rho_k a\nabla u)+\rho_k cu,\quad\mbox{for}\,\,k=1,2,\ldots,m.
\end{equation}
The family of functions $\{\rho_k(\mathbf{x})\}_{k=1,2,\ldots,m}$ is subordinate to an overlapping decomposition of $\Omega$ and conforms a smooth partition of unity. The construction of such a partition is discussed in the sequel.

Let $\{\Omega_k^{\ast}\}_{k=1,2,\ldots,m}$ form a non-overlapping decomposition of $\Omega$ into $m$ subdomains. Such a decomposition fulfills the conditions $\overline{\Omega}=\bigcup_{k=1}^m\overline{\Omega}^{\ast}_k$ and $\Omega_{k}^{\ast}\cap\Omega_{l}^{\ast}=\emptyset$, for $k\neq l$. In turn, each $\Omega_{k}^{\ast}\subset\Omega$ is considered to be an open disconnected set involving $m_k$ connected components, i.e.,
\begin{equation*}
\textstyle\Omega_k^{\ast}=\bigcup_{l=1}^{m_{k}}\Omega^{\ast}_{kl},\quad\mathrm{for}\,\,k=1,2,\ldots,m.
\end{equation*}
Such components are pairwise disjoint (that is, $\Omega_{ki}^{\ast}\cap\Omega_{kj}^{\ast}=\emptyset$, for $i\neq j$) and typically chosen to be shape regular of diameter $h_0$. For instance, the components $\Omega_{kl}^{\ast}$ may correspond to the elements in a coarse partition $\mathcal{T}_{h_0}$ of $\Omega$ with mesh size $h_0$. Let $\Omega_{kl}$ be the extension of $\Omega_{kl}^{\ast}$ obtained by suitably translating its internal boundaries, $\partial\Omega_{kl}^{\ast}\cap\Omega$, within a distance $\beta h_0$ in $\Omega$. The parameter $\beta>0$ is usually referred to as the overlapping factor and its value must be chosen in such a way that the extended components are also pairwise disjoint (i.e., $\Omega_{ki}\cap\Omega_{kj}=\emptyset$, for $i\neq j$). The distance $\xi=2\beta h_0$ is called the overlapping size. If we denote by $\Omega_k\subset\Omega$ the open disconnected set defined as
\begin{equation}\label{disjoint:components}
\textstyle\Omega_k=\bigcup_{l=1}^{m_{k}}\Omega_{kl},\quad\mathrm{for}\,\,k=1,2,\ldots,m,
\end{equation}
then the collection $\{\Omega_k\}_{k=1,2,\ldots,m}$ forms an overlapping decomposition of $\Omega$ into $m$ subdomains, that is, $\Omega=\bigcup_{k=1}^m\Omega_k$.

Subordinate to this overlapping covering of $\Omega$, we construct a smooth partition of unity consisting of a family of $m$ non-negative and $\mathcal{C}^{\infty}(\Omega)$ functions $\{\rho_k(\mathbf{x})\}_{k=1,2,\ldots,m}$. Each function $\rho_{k}:\overline{\Omega}\rightarrow[0,1]$ is chosen to be
\begin{equation*}
\rho_{k}(\mathbf{x})=
\begin{cases}
\,0,&\mathrm{if}\,\,\mathbf{x}\in\overline{\Omega}\setminus\overline{\Omega}_{k},\\
\,h_{k}(\mathbf{x}),&\mathrm{if}\,\,\mathbf{x}\in\bigcup_{l=1;\,l\neq k}^{m}\,(\overline{\Omega}_{k}\cap\overline{\Omega}_l),\\
\,1,&\mathrm{if}\,\,\mathbf{x}\in\overline{\Omega}_{k}\setminus\bigcup_{l=1;\,l\neq k}^{m}\,(\overline{\Omega}_{k}\cap\overline{\Omega}_l),
\end{cases}
\end{equation*}
where $h_k(\mathbf{x})$ is $\mathcal{C}^{\infty}(\Omega)$ and satisfies the conditions
\begin{equation*}
0\leq h_k(\mathbf{x})\leq1,\qquad\sum_{k=1}^mh_k(\mathbf{x})=1,
\end{equation*}
for any $\mathbf{x}\in\bigcup_{l=1;\,l\neq k}^{m}\,(\overline{\Omega}_{k}\cap\overline{\Omega}_l)$. The family of functions $\{\rho_k(\mathbf{x})\}_{k=1,2,\ldots,m}$ is such that
\begin{equation}\label{partition:of:unity:conditions}
\mathrm{supp}(\rho_k(\mathbf{x}))\subset\overline{\Omega}_k,\qquad
0\leq\rho_k(\mathbf{x})\leq1,\qquad
\sum_{k=1}^m\rho_k(\mathbf{x})=1,
\end{equation}
for any $\mathbf{x}\in\overline{\Omega}$. For numerical purposes, $h_k(\mathbf{x})$ may not necessarily be $\mathcal{C}^{\infty}(\Omega)$, but only a continuous and piecewise smooth function (cf. \cite{mat:pol:rus:wan:98}); this fact will be illustrated in the numerical experiments provided below.

In component-wise splittings, the number $m$ of split terms is equal to the dimension $d$ of the spatial domain $\Omega$. In particular, the alternating direction splitting corresponding to the two-dimensional problem (\ref{continuous:problem}) consists of two split terms given by (\ref{adi:splitting}). However, for domain decomposition splittings as that considered in (\ref{dd:splitting}), $m$ is not necessarily equal to $d$. As a result, such splittings may require a generalization of the two-stage procedure (\ref{adi:scheme}) to an arbitrary number $m$ of stages. This general formulation was first proposed in \cite{dou:gun:64} in the context of ADI methods. In our case, it is given, for $n=0,1,\ldots,N_T$ and $k=1,2,\ldots,m$, by
\begin{equation}\label{dd:scheme}
\begin{aligned}
\frac{W_h^{n,k}-W_h^n}{\tau}+\sum_{i=1}^{k}A_{ih}(\theta W_h^{n,i}+(1-\theta)W_h^n)+\sum_{i=k+1}^{m}A_{ih}W_h^n&=F_h^{n+\theta},\\
W_h^{n+1}&=W_h^{n,m},
\end{aligned}
\end{equation}
where $W_h^0=\mathcal{P}_hu_0$. Note that the domain decomposition discrete solution is denoted by $W_h^n$. Here, $\{A_{kh}\}_{k=1,2,\ldots,m}$ are suitable discretizations of the split operators $\{A_{k}\}_{k=1,2,\ldots,m}$ introduced in (\ref{dd:splitting}); such discretizations are constructed to satisfy $A_h=A_{1h}+A_{2h}+\ldots+A_{mh}$. For convenience, the algorithm (\ref{dd:scheme}) is rewritten in the equivalent form
\begin{subequations}\label{generalized:splitting}
\begin{align}
(I+\theta\tau A_{1h})\,W_h^{n,1}&=\left(I-(1-\theta)\tau A_{1h}-\tau\sum_{i=2}^mA_{ih}\right)W_h^n+\tau F_h^{n+\theta},\label{generalized:splitting:a}\\[0.5ex]
(I+\theta\tau A_{kh})\,W_h^{n,k}&=W_h^{n,k-1}+\theta\tau A_{kh}W_h^n,\quad\mbox{for\,\,}k=2,3,\ldots,m,\label{generalized:splitting:b}\\[2.2ex]
W_h^{n+1}&=W_h^{n,m}.\label{generalized:splitting:c}
\end{align}
\end{subequations}
Thus, the linear system to be solved at the $k$-th internal stage is of the form
\begin{equation}\label{global:system}
(I+\theta\tau{A}_{kh})W_h^{n,k}=Q_h^{n,k},\quad\mbox{for\,\,}k=1,2,\ldots,m,
\end{equation}
where $Q_h^{n,k}$ stands for the corresponding right-hand side. As stated in (\ref{partition:of:unity:conditions}), the function $\rho_{k}(\mathbf{x})$ has compact support on $\overline{\Omega}_{k}$ and, by construction, the entries of $A_{kh}$ corresponding to the nodes that lie outside of this subdomain are zero. In addition, since $\Omega_k$ involves $m_k$ disjoint connected components $\Omega_{kl}$ due to (\ref{disjoint:components}), the linear system (\ref{global:system}) is indeed a collection of $m_k$ uncoupled subsystems of the form
\begin{equation}\label{local:system}
(I_{\Omega_{kl}}+\theta\tau{A}_{\Omega_{kl}h})\mathcal{R}_{\Omega_{kl}}W_h^{n,k}=\mathcal{R}_{\Omega_{kl}}Q_h^{n,k},\quad\mbox{for\,\,}k=1,2,\ldots,m,
\end{equation}
where $I_{\Omega_{kl}}=\mathcal{R}_{\Omega_{kl}}I\mathcal{R}_{\Omega_{kl}}^T$ and ${A}_{\Omega_{kl}h}=\mathcal{R}_{\Omega_{kl}}A_{kh}\mathcal{R}_{\Omega_{kl}}^T$. The rectangular matrices $\mathcal{R}_{\Omega_{kl}}$ and $\mathcal{R}_{\Omega_{kl}}^T$ are usually called restriction and extension matrices, respectively, and represent a type of domain decomposition preconditioners (cf. \cite{mat:08}). To obtain such matrices, we first order the nodes of $\mathcal{T}_h$ in $\Omega_{kl}$ in some local ordering. Let us denote by $N_{\mathcal{S}}$ and $N$, respectively, the number of nodes of $\mathcal{T}_h$ in $\Omega_{kl}$ and the total number of nodes of $\mathcal{T}_h$. Then, $\mathcal{R}_{\Omega_{kl}}$ is an $N_{\mathcal{S}}\times N$ matrix defined as
\begin{equation*}
(\mathcal{R}_{\Omega_{kl}})_{ij}=\begin{cases}
\,1,&\mathrm{if}\,\,{index}(\Omega_{kl},i)=j,\\
\,0,&\mathrm{if}\,\,{index}(\Omega_{kl},i)\neq j,
\end{cases}
\end{equation*}
where the index function \emph{index}$(\Omega_{kl},i)$ obtains the global index of the $i$-th local node of $\mathcal{T}_h$ in $\Omega_{kl}$, for $i=1,2,\ldots,N_{\mathcal{S}}$. In other words, this matrix maps a vector in $\mathds{R}^{N}$ of nodal values in $\mathcal{T}_h$ into a vector in $\mathds{R}^{N_{\mathcal{S}}}$ of nodal values corresponding to the nodes of $\mathcal{T}_h$ in $\Omega_{kl}$ (in the local ordering). Likewise, the $N\times N_{\mathcal{S}}$ transpose matrix $\mathcal{R}_{\Omega_{kl}}^{T}$ maps a vector in $\mathds{R}^{N_{\mathcal{S}}}$ into a vector in $\mathds{R}^{N}$, inserting zero entries for the global indices which do not belong to $\Omega_{kl}$. Finally, given a global matrix $M_h\in\mathds{R}^{N\times N}$, the submatrix $M_{\Omega_{kl}h}\in\mathds{R}^{N_{\mathcal{S}}\times N_{\mathcal{S}}}$ associated with the nodes in $\Omega_{kl}$ may be obtained from the restriction and extension mappings as $M_{\Omega_{kl}h}=\mathcal{R}_{\Omega_{kl}}M_h\mathcal{R}_{\Omega_{kl}}^{T}$.

From a computational viewpoint, this domain decomposition splitting procedure has three main parameters to adjust its efficiency, namely: the number $m$ of subdomains, the number $m_k$ of disjoint connected components inside a subdomain, and the overlapping size $\xi$. In order to achieve an efficient performance of the algorithm, $m$ should be chosen as small as possible to minimize the number of sequential steps (i.e., internal stages) in (\ref{local:system}), $m_k$ should be chosen as large as possible to maximize the number of parallel components, and $\xi$ should be chosen as small as possible to minimize the number of unknowns within the overlapping regions. In addition, to ensure that the loads assigned to each processor are balanced, there should be approximately the same number of disjoint connected components in each subdomain and, moreover, each such component should have approximately the same diameter. For instance, suppose that we have $q$ available processors for parallel computing. Then, if $m_k$ is approximately the same in each subdomain $\Omega_k$ and also a multiple of $q$, each subdomain can be partitioned into $q$ groups of $\frac{m_k}{q}$ components and each group assigned to one of the processors. Remarkably, unlike most classical overlapping domain decomposition algorithms (cf. \cite{qua:val:99}), the solution to (\ref{local:system}) does not require any Schwarz iteration procedure, since the internal stages in the algorithm are solved sequentially (i.e., interface conditions need not be imposed on subdomains during the solution process).

\section{The splitting error}\label{sec:splitting:error}

As noted in the introduction, the time-splitting methods presented above can be interpreted as perturbations of some unsplit implicit schemes, such as the backward Euler and Crank--Nicolson methods. Recall that the scheme (\ref{adi:scheme}) considers the alternating direction splitting (\ref{adi:splitting}), while the algorithm (\ref{dd:scheme}) makes use of the domain decomposition splitting (\ref{dd:splitting}). In essence, the former method is a particular instance of the latter for $m=2$. In this section, we reformulate (\ref{dd:scheme}) as a classical implicit scheme (backward Euler, for $\theta=1$, and Crank--Nicolson, for $\theta=\frac{1}{2}$) plus a perturbation term, which will be referred to as the splitting error. As we will see, this splitting error is $\mathcal{O}(\tau^2)$, i.e., one order higher than the truncation error associated with the backward Euler method and of the same order as that associated with the Crank--Nicolson method. As a result, one would expect to obtain the same asymptotic rate of convergence for both the time-splitting and its underlying implicit method. In practice, however, the size of the splitting error can be much larger than that of the truncation error corresponding to the underlying scheme, thus degrading the efficiency of the resulting algorithm. In the sequel, we derive a general expression for the splitting error and present a strategy to reduce its actual size.

Let us consider the formulation given by (\ref{generalized:splitting}). In order to eliminate the intermediate values $W_h^{n,1},W_h^{n,2},\ldots,W_h^{n,m-1}$, we multiply the equations (\ref{generalized:splitting:b}) by $(I+\theta\tau A_{1h})(I+\theta\tau A_{2h})\cdots(I+\theta\tau A_{(k-1)h})$, sum on $k$ and add the equation (\ref{generalized:splitting:a}) to obtain, for $n=0,1,\ldots,N_T$,
\begin{equation}\label{perturbed:implicit}
\frac{W_h^{n+1}-W_h^n}{\tau}+A_h(\theta W_h^{n+1}+(1-\theta)W_h^n)+B_h(W_h^{n+1}-W_h^n)=F_h^{n+\theta},
\end{equation}
where
\begin{align*}
B_h=&\,\,\theta^2\tau\sum_{1\leq i_1<i_2\leq m}A_{i_1h}A_{i_2h}+\theta^3\tau^2\sum_{1\leq i_1<i_2<i_3\leq m}A_{i_1h}A_{i_2h}A_{i_3h}\\
&\,\,+\ldots+\theta^m\tau^{m-1}A_{1h}A_{2h}\cdots A_{mh}.
\end{align*}
Note that, depending on the value of $\theta$, (\ref{perturbed:implicit}) is, in fact, a perturbed backward Euler or Crank--Nicolson method (as compared with (\ref{be:cn})). The perturbation term $B_h(W_h^{n+1}-W_h^n)$ is the aforementioned splitting error and satisfies
\begin{equation*}
B_h(W_h^{n+1}-W_h^n)=\tau B_h\left(\frac{W_h^{n+1}-W_h^n}{\tau}\right)=\mathcal{O}\left(\tau^2\sum_{i,j}\left|\frac{\partial}{\partial t}A_iA_ju\right|\right)=\mathcal{O}(\tau^2),
\end{equation*}
where the split operators $A_k$ are given by either (\ref{adi:splitting}) or (\ref{dd:splitting}), provided that $u$, $c$ and the entries of $a$ are sufficiently smooth.

In \cite{dou:kim:01}, Douglas and Kim noted that, if we could add $B_h(W_h^{n+1}-W_h^n)$ to the right-hand side of (\ref{perturbed:implicit}), then the perturbation term would be cancelled. However, since $W_h^{n+1}$ is not known at time $t_n$, this modification cannot be done. Taking into account that the best available approximation to the difference $W_h^{n+1}-W_h^n$ at time $t_n$ is $W_h^{n}-W_h^{n-1}$, the algorithm (\ref{generalized:splitting}) can be modified, for $n=1,2,\ldots,N_T$, to obtain
\begin{equation}\label{improved:splitting}
\begin{aligned}
\widehat{F}_h^{n+\theta}&=F_h^{n+\theta}+B_h(Z_h^{n}-Z_h^{n-1}),\\[0.5ex]
(I+\theta\tau A_{1h})Z_h^{n,1}&=\left(I-(1-\theta)\tau A_{1h}-\tau\sum_{i=2}^mA_{ih}\right)Z_h^n+\tau \widehat{F}_h^{n+\theta}\\[0.5ex]
(I+\theta\tau A_{kh})Z_h^{n,k}&=Z_h^{n,k-1}+\theta\tau A_{kh}Z_h^n,\qquad\mbox{for\,\,}k=2,3,\ldots,m,\\[2.2ex]
Z_h^{n+1}&=Z_h^{n,m},
\end{aligned}
\end{equation}
where $Z_h^0=\mathcal{P}_hu_0$ and $Z_h^1$ must be suitably approximated. In this case, eliminating the intermediate values $Z_h^{n,1},Z_h^{n,2},\ldots,Z_h^{n,m-1}$, or referring to (\ref{perturbed:implicit}), we get, for $n=1,2,\ldots,N_T$,
\begin{equation*}
\frac{Z_h^{n+1}-Z_h^n}{\tau}+A_h(\theta Z_h^{n+1}+(1-\theta)Z_h^n)+B_h(Z_h^{n+1}-2Z_h^n+Z_h^{n-1})=F_h^{n+\theta}.
\end{equation*}
Now, the perturbation term $B_h(Z_h^{n+1}-2Z_h^n+Z_h^{n-1})$ satisfies
\begin{align*}
B_h(Z_h^{n+1}-2Z_h^n+Z_h^{n-1})&=\tau^2 B_h\left(\frac{Z_h^{n+1}-2Z_h^n+Z_h^{n-1}}{\tau^2}\right)\\
&=\mathcal{O}\left(\tau^3\sum_{i,j}\left|\frac{\partial^2}{\partial t^2}A_iA_ju\right|\right)=\mathcal{O}(\tau^3),
\end{align*}
provided that $u$, $c$ and the entries of $a$ are sufficiently smooth. In other words, the corrected right-hand side $\widehat{F}_h^{n+\theta}$ introduced in (\ref{improved:splitting}) yields a reduction in the order of the splitting error from $\mathcal{O}(\tau^2)$ to $\mathcal{O}(\tau^3)$. As a result, the actual discretization error of the modified algorithm (\ref{improved:splitting}) is of the same size as that associated with its underlying unsplit method. Remarkably, if the domain decomposition splitting (\ref{dd:splitting}) is used to obtain $A_k$, the operator $B_h$ will be non-zero only within the overlapping regions. Thus, the newly added term $B_h(Z_h^{n}-Z_h^{n-1})$ will have a very low extra computational cost.

The convergence analysis for the improved accuracy method (\ref{improved:splitting}) must take into account the splitting strategy under consideration. The case in which the split operators $A_k$ are obtained via the alternating direction splitting (\ref{adi:splitting}) is discussed in \cite{arb:hua:yan:07} and \cite{dou:kim:01}, and error estimates of order $\mathcal{O}(\tau^{3-2\theta}+h^s)$ are derived. In the latter work, the discrete split operators are assumed to commute pairwise, while the former analyzes the non-commuting case. As for the domain decomposition splitting (\ref{dd:splitting}), a convergence analysis of this kind is lacking so far. A similar analysis for the unimproved algorithm (\ref{generalized:splitting}), considering a finite difference spatial discretization, is carried out in \cite{mat:pol:rus:wan:98}. The derivation of error estimates for the newly proposed method (\ref{improved:splitting}) is a topic of current research.

\section{Numerical experiments}\label{sec:experiments}

The numerical experiments performed in this section are partly inspired by those contained in \cite{arb:hua:yan:07,dou:kim:01}. Let us consider the parabolic initial-boundary value problem (\ref{continuous:problem}), where $\Omega=(0,1)\times(0,1),$ $T=1$ and $c(\mathbf{x})\equiv 0$. Given a certain diffusion term $a(\mathbf{x})$ (to be specified below), the functions $f(\mathbf{x},t)$ and $u_0(\mathbf{x})$ are defined in such a way that the exact solution is
\begin{equation*}
u(x,y,t)=\sin(2\pi t)\sin(2\pi x)\sin(2\pi y).
\end{equation*}
Throughout this section, a finite difference spatial discretization is combined with several second-order time integrators (i.e., we fix $\theta=\frac{1}{2}$), considering the discretization parameters $h=\tau=\frac{1}{M}$. In particular, we compare the accuracy of the Crank--Nicolson (CN) method (\ref{be:cn}) to that of the Douglas--Gunn scheme (\ref{generalized:splitting}) and the Douglas--Kim method (\ref{improved:splitting}), combined with either the ADI splitting (\ref{adi:splitting}) (DG$_{\mathrm{ADI}}$ and DK$_{\mathrm{ADI}}$, respectively) or the domain decomposition splitting (\ref{dd:splitting}) (DG$_{\mathrm{DD}}$ and DK$_{\mathrm{DD}}$, respectively). For the Douglas--Kim procedures, DK$_{\mathrm{ADI}}$ and DK$_{\mathrm{DD}}$, the approximation of the exact solution at time $t_1=\tau$ is obtained by running a single step of the Crank--Nicolson scheme. The computed errors are measured in the discrete $\ell^{\infty}(\ell^2)$-norm, that is,
\begin{equation*}
\|e\|_{\ell^{\infty}(\ell^2)}=\max_{1\leq n\leq N_T}\left(\sum_{i=1}^{M-1}\sum_{j=1}^{M-1}(e_{i,j}^n)^2h^2\right)^{1/2},
\end{equation*}
where $e_{i,j}^n$ denotes the error at the $(i,j)$-th node and time $t_n$, obtained as the difference between the exact solution at $(ih,jh,t_n)$ and the corresponding numerical solution.

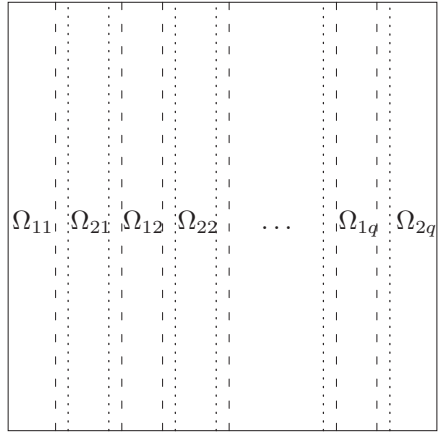
\begin{figure}[t]
\begin{center}
\unitlength=0.057cm
\begin{picture}(100,100)
\thinlines
\path(0,0)(100,0)(100,100)(0,100)(0,0)
\dashline[-10]{2}(11,0)(11,100)
\dottedline{2}(14,0)(14,100)
\dottedline{2}(23.5,0)(23.5,100)
\dashline[-10]{2}(26.5,0)(26.5,100)
\dashline[-10]{2}(36,0)(36,100)
\dottedline{2}(39,0)(39,100)
\dottedline{2}(48.5,0)(48.5,100)
\dashline[-10]{2}(51.5,0)(51.5,100)
\dottedline{2}(73.5,0)(73.5,100)
\dashline[-10]{2}(76.5,0)(76.5,100)
\dashline[-10]{2}(86,0)(86,100)
\dottedline{2}(89,0)(89,100)
\put(1,47){\footnotesize{$\Omega_{11}$}}
\put(14.5,47){\footnotesize{$\Omega_{21}$}}
\put(26.7,47){\footnotesize{$\Omega_{12}$}}
\put(39.5,47){\footnotesize{$\Omega_{22}$}}
\put(59,47){\footnotesize{$\ldots$}}
\put(77,47){\footnotesize{$\Omega_{1q}$}}
\put(90.5,47){\footnotesize{$\Omega_{2q}$}}
\end{picture}
\end{center}
\caption{Decomposition of $\Omega$ into two subdomains, $\Omega_1$ and $\Omega_2$, each consisting of $q$ disjoint connected components, $\{\Omega_{1l}\}_{l=1,2,\ldots,q}$ and $\{\Omega_{2l}\}_{l=1,2,\ldots,q}$. The internal boundaries of such components are represented by dotted and dashed lines, respectively.}
\label{fig:1D:DD:splitting}
\end{figure}

Unless otherwise stated, both the DG$_{\mathrm{DD}}$ and the DK$_{\mathrm{DD}}$ methods consider two overlapping subdomains, each consisting of four non-overlapping connected components. More precisely, let $I=(0,1)$ and define the intervals
\begin{align*}
I_1&=\textstyle\left(0,\frac{1}{8}+\frac{\xi}{2}\right)\cup\left(\frac{1}{4}-\frac{\xi}{2},\frac{3}{8}+\frac{\xi}{2}\right)\cup\left(\frac{1}{2}-\frac{\xi}{2},\frac{5}{8}+\frac{\xi}{2}\right)
\cup\left(\frac{3}{4}-\frac{\xi}{2},\frac{7}{8}+\frac{\xi}{2}\right),\\[0.5ex]
I_2&=\textstyle\left(\frac{1}{8}-\frac{\xi}{2},\frac{1}{4}+\frac{\xi}{2}\right)\cup\left(\frac{3}{8}-\frac{\xi}{2},\frac{1}{2}+\frac{\xi}{2}\right)\cup\left(\frac{5}{8}-\frac{\xi}{2},\frac{3}{4}+\frac{\xi}{2}\right)
\cup\left(\frac{7}{8}-\frac{\xi}{2},1\right),
\end{align*}
where $\xi$ denotes the overlapping size (to be specified below). Then, we can define the subdomains $\Omega_1=I_1\times I$ and $\Omega_2=I_2\times I$. Figure \ref{fig:1D:DD:splitting} shows a domain decomposition of this kind on the unit square, for an arbitrary number $q$ of connected components per subdomain.


According to \cite{mat:pol:rus:wan:98}, we construct a piecewise smooth partition of unity as follows. The closure of subdomain $\Omega_k$ can be expressed as the union of four closed disjoint rectangles given by $R_k^i=[a_k^i,b_k^i]\times[0,1]$, for $i=1,2,3,4$ and $k=1,2$. Then, for each rectangle $R_k^i$, we define a function of the form
\begin{equation*}
w_{k}^i(x,y)=
\begin{cases}
\,\sin\left(\displaystyle\frac{\pi(x-a_k^i)}{b_k^i-a_k^i}\right),&\mbox{if }(x,y)\in R_k^i,\\
\,0,&\mbox{if }(x,y)\in\overline{\Omega}\setminus R_k^i.
\end{cases}
\end{equation*}
The partition of unity functions are thus obtained as
\begin{equation}\label{piecewise:smooth:partition:unity}
\rho_k(x,y)=\frac{\displaystyle\sum_{i=1}^{4}w_k^i(x,y)}{\displaystyle\sum_{l=1}^{2}\displaystyle\sum_{i=1}^{4}w_l^i(x,y)},\quad\mbox{for } k=1,2.
\end{equation}

Table \ref{table:N:varies} shows the global errors and mean rates of convergence obtained for different values of the discretization parameter $M$, with the aforementioned methods. In this experiment, the diffusion coefficient is $a(x,y)=1$, and the overlapping size is $\xi=\frac{1}{8}$. We observe that the global errors corresponding to the DG$_{\mathrm{ADI}}$ and the DG$_{\mathrm{DD}}$ methods are about 10 and 14 times larger, respectively, than the CN errors. The DK$_{\mathrm{ADI}}$ errors are larger but comparable to the CN errors, since the splitting error for the former method is much smaller than that for the DG$_{\mathrm{ADI}}$ scheme. In this case, the mean rate of convergence is greater than 2, because the splitting error is being removed at a rate of $\mathcal{O}(\tau^3)$. Accordingly, the DK$_{\mathrm{DD}}$ errors are much smaller than the DG$_{\mathrm{DD}}$ errors, and very similar to those obtained for the CN scheme. In fact, due to the observed higher rate of convergence, the DK$_{\mathrm{DD}}$ errors are indeed smaller than the CN errors, for $M=160$ and $M=320$.


In Table \ref{table:a:varies}, we fix the values $M=160$ and $\xi=\frac{1}{8}$, and compute the global errors corresponding to the previous schemes, for five different choices of the diffusion coefficient $a(x,y)$. Besides the one tested in the preceding table, $a_1(x,y)=1$, we also consider
\begin{align*}
a_2(x,y)&=\displaystyle\frac{1}{2+\cos(3\pi x)\cos(2\pi y)},\\[1.75ex]
a_3(x,y)&=
\begin{cases}
\,1+0.5 \sin(5\pi x)+y^3,&\mbox{if }x\leq 0.5,\\
\,\dfrac{1.5}{1+(x-0.5)^2}+y^3,&\mbox{otherwise},
\end{cases}\\[1.75ex]
a_4(x,y)&=\left(\begin{array}{cc}a_2(x,y) &0\\0&a_3(x,y)\end{array}\right),\\[1.75ex]
a_5(x,y)&=\left(\begin{array}{cc}a_2(x,y) &1/4\\1/4&a_2(x,y)\end{array}\right).
\end{align*}
If we take $a(x,y)=a_i(x,y)$, for $i=1,2,3,4$, the computed DG$_{\mathrm{ADI}}$ errors are much larger than the corresponding CN errors. In turn, the errors associated with the DK$_{\mathrm{ADI}}$ method are comparable to, and only slightly larger than, those obtained with the CN scheme. On the other hand, if the diffusion coefficient is chosen to be $a_5(x,y)$, none of these two ADI methods can be applied, due to the presence of mixed derivative terms. As for the domain decomposition splitting techniques, in all five cases, the DK$_{\mathrm{DD}}$ scheme is shown to be comparable to the CN scheme and superior to the DG$_{\mathrm{DD}}$ method.


Table \ref{table:d:varies} shows the effect of the overlapping size $\xi$ on the accuracy of the implemented methods. In particular, we consider $a(x,y)=a_2(x,y)$ and $M=160$, and compute the global errors corresponding to the CN, the DG$_{\mathrm{DD}}$ and the DK$_{\mathrm{DD}}$ schemes, for different values of $\xi$. The observed DG$_{\mathrm{DD}}$ errors increase as $\xi$ becomes smaller, due to the presence of negative powers of $\xi$ in the splitting error (cf. \cite{mat:pol:rus:wan:98}). This undesirable behaviour is avoided by the DK$_{\mathrm{DD}}$ scheme, whose global errors get slightly smaller as $\xi$ decreases.

Finally, we study the influence of the number of disjoint connected components per subdomain on the accuracy of the domain decomposition procedures. For that purpose, we consider the same number $q$ of components in $\Omega_1$ and $\Omega_2$; i.e., recalling the notation introduced in \S\ref{sec:split}, $m_1=m_2=q$. Table \ref{table:q:varies} shows the global errors corresponding to the CN, the DG$_{\mathrm{DD}}$ and the DK$_{\mathrm{DD}}$ schemes, for different values of $q$, when $a(x,y)=a_2(x,y)$, $M=160$ and $\xi=\frac{1}{16}$. Note that, if $q=2$ or $q=8$, the partition of unity functions (\ref{piecewise:smooth:partition:unity}) must be modified accordingly. As in the preceding experiments, the DG$_{\mathrm{DD}}$ errors are larger than the CN errors; moreover, they even grow as $q$ increases, which represents a serious drawback from a computational viewpoint. By contrast, the errors for the DK$_{\mathrm{DD}}$ scheme are comparable to, and even smaller than, those associated with the CN method. In this case, such errors decrease as $q$ increases, thus permitting to exploit the parallel capabilities of the new algorithms.



\section{Conclusions}\label{sec:conclusions}

A family of improved domain decomposition splitting methods for solving parabolic equations has been presented. These methods can be formulated as unsplit implicit schemes (such as backward Euler or Crank--Nicolson) perturbed by a term that is $\mathcal{O}(\tau^3)$ in the time step $\tau$. This perturbation term, also known as splitting error, is one order higher than its $\mathcal{O}(\tau^2)$ analogue stemming from the unimproved formulation. Such an increase in the order of the splitting error is achieved by adding a correction term to the right-hand side of the original scheme. This technique was proposed in \cite{dou:kim:01} in the context of ADI methods, and has been extended here to domain decomposition-based splittings for the first time. The newly derived methods can be implemented with a very low additional cost, since the correction term only requires certain computations within the overlapping regions. In addition, the resulting subdomain problems can be solved in parallel with no need for Schwarz-type iteration procedures. Numerical results show that the proposed modification leads to a significant reduction in the splitting error. Moreover, the global discretization error is comparable to that associated with the Crank--Nicolson method, independently of the number of connected components per subdomain and the overlapping size.

%


\clearpage

\begin{table}[t]
\begin{center}
\renewcommand{\arraystretch}{1.2}
\centering{
\begin{tabular}{l|c|c|c|c|c} \hline
Method & $M=40$ &$M=80$ &$M=160$ &$M=320$ & Rate \\
\hhline{|------|}
CN & 1.029e-03 & 2.571e-04 & 6.426e-05 & 1.606e-05 & 2.000 \\
DG$_{\mathrm{ADI}}$ & 9.780e-03 & 2.457e-03 & 6.244e-04 & 1.658e-04 & 1.961 \\
DK$_{\mathrm{ADI}}$ & 2.558e-03 & 4.560e-04 & 9.920e-05 & 3.121e-05 & 2.119 \\
DG$_{\mathrm{DD}}$ & 1.444e-02 & 3.026e-03 & 8.488e-04 & 2.252e-04 & 2.001 \\
DK$_{\mathrm{DD}}$ & 2.180e-03 & 2.933e-04 & 6.079e-05 & 1.494e-05 & 2.396 \\
\hline
\end{tabular}}
\end{center}
\caption{Global errors and mean rates of convergence for different values of the discretization parameter $M$, with $a(x,y)=1$ and $\xi=\frac{1}{8}$.}
\label{table:N:varies}
\end{table}


\begin{table}[t]
\begin{center}
\renewcommand{\arraystretch}{1.2}
\centering{
\begin{tabular}{l|c|c|c|c|c} \hline
Method  & $a=a_1$ &$a=a_2$ &$a=a_3$ &$a=a_4$ &$a=a_5$ \\
\hhline{|------|}
CN & 6.426e-05 & 6.179e-05 & 7.456e-05 & 6.160e-05 & 9.339e-05 \\
DG$_{\mathrm{ADI}}$ & 6.244e-04 & 3.541e-04 & 9.391e-04 & 4.863e-04 & -- \\
DK$_{\mathrm{ADI}}$ & 9.920e-05 & 7.732e-05 & 1.116e-04 & 8.308e-05 & -- \\
DG$_{\mathrm{DD}}$ & 8.488e-04 & 4.427e-04 & 1.464e-03 & 5.313e-04 & 5.333e-04 \\
DK$_{\mathrm{DD}}$ & 6.079e-05 & 5.905e-05 & 8.327e-05 & 6.481e-05 & 9.593e-05 \\
\hline
\end{tabular}}
\end{center}
\caption{Global errors for different choices of the diffusion coefficient $a(x,y)$, with $M=160$ and $\xi=\frac{1}{8}$.}
\label{table:a:varies}
\end{table}


\begin{table}[t]
\begin{center}
\renewcommand{\arraystretch}{1.2}
\centering{
\begin{tabular}{l|c|c|c} \hline
Method  & $\xi=\frac{1}{8}$ &$\xi=\frac{1}{16}$ &$\xi=\frac{1}{32}$ \\
\hhline{|----|}
CN & 6.179e-05 & 6.179e-05 & 6.179e-05 \\
DG$_{\mathrm{DD}}$ & 4.427e-04& 8.148e-04 & 1.747e-03 \\
DK$_{\mathrm{DD}}$ & 5.905e-05 & 4.600e-05 & 4.418e-05 \\
\hline
\end{tabular}}
\end{center}
\caption{Global errors for different values of the overlapping size $\xi$, with $a(x,y)=a_2(x,y)$ and $M=160$.}
\label{table:d:varies}
\end{table}


\begin{table}[t]
\begin{center}
\renewcommand{\arraystretch}{1.2}
\centering{
\begin{tabular}{l|c|c|c} \hline
Method  & $q=2$ &$q=4$ &$q=8$ \\
\hhline{|----|}
CN & 6.179e-05 & 6.179e-05 & 6.179e-05 \\
DG$_{\mathrm{DD}}$ & 4.407e-04& 8.148e-04 & 1.475e-03 \\
DK$_{\mathrm{DD}}$ & 5.717e-05 & 4.600e-05 & 3.742e-05 \\
\hline
\end{tabular}}
\end{center}
\caption{Global errors for different number $q$ of connected components per subdomain, with $a(x,y)=a_2(x,y)$, $M=160$ and $\xi=\frac{1}{16}$.}
\label{table:q:varies}
\end{table}


\end{document}